\def\hldest#1#2#3{}



\font\eightpt=cmr8
\font\ninept=cmr9


\font\ninerm=cmr9     \font\eightrm=cmr8   \font\sixrm=cmr6      
\font\ninei=cmmi9     \font\eighti=cmmi8   \font\sixi=cmmi6      
\font\ninesy=cmsy9    \font\eightsy=cmsy8  \font\sixsy=cmsy6     
\font\ninebf=cmbx9    \font\eightbf=cmbx8  \font\sixbf=cmbx6     
\font\ninett=cmtt9    \font\eighttt=cmtt8                        
\font\nineit=cmti9    \font\eightit=cmti8     
\font\ninesl=cmsl9    \font\eightsl=cmsl8                        

\font\tensc=cmcsc10   \font\ninesc=cmcsc9  \font\eightsc=cmcsc8  

\font\eightssq=cmssq8  \font\eightssqi=cmssqi8  


\font\tenssbx=cmssbx10 

  \font\twelvebf=cmbx12
  
\def\sc{\tensc}  \def\mc{\ninerm}


\newskip\ttglue
\def\tenpoint{\def\rm{\fam0\tenrm}%
  \textfont0=\tenrm \scriptfont0=\sevenrm \scriptscriptfont0=\fiverm
  \textfont1=\teni  \scriptfont1=\seveni  \scriptscriptfont1=\fivei
  \textfont2=\tensy \scriptfont2=\sevensy \scriptscriptfont2=\fivesy
  \textfont3=\tenex \scriptfont3=\tenex   \scriptscriptfont3=\tenex
  \textfont\itfam=\tenit  \def\it{\fam\itfam\tenit}%
  \textfont\slfam=\tensl  \def\sl{\fam\slfam\tensl}%
  \textfont\ttfam=\tentt  \def\tt{\fam\ttfam\tentt}%
  \textfont\bffam=\tenbf  \scriptfont\bffam=\sevenbf
   \scriptscriptfont\bffam=\fivebf \def\bf{\fam\bffam\tenbf}%
  \tt \ttglue=.5em plus.25em minus.15em
  \normalbaselineskip=12pt
  \setbox\strutbox=\hbox{\vrule height8.5pt depth3.5pt width0pt}%
  \let\sc=\tensc \let\mc=\ninerm  
  \def\cyr{\tencyr\cyracc}\def\cyri{\tencyri\cyracc}
  \let\big=\tenbig  \normalbaselines\rm}

\def\ninepoint{\def\rm{\fam0\ninerm}%
\textfont0=\ninerm  \scriptfont0=\sixrm  \scriptscriptfont0=\fiverm
\textfont1=\ninei   \scriptfont1=\sixi   \scriptscriptfont1=\fivei
\textfont2=\ninesy  \scriptfont2=\sixsy  \scriptscriptfont2=\fivesy
\textfont3=\tenex   \scriptfont3=\tenex  \scriptscriptfont3=\tenex
\textfont\itfam=\nineit  \def\it{\fam\itfam\nineit}%
\textfont\slfam=\ninesl  \def\sl{\fam\slfam\ninesl}%
\textfont\ttfam=\ninett  \def\tt{\fam\ttfam\ninett}%
\textfont\bffam=\ninebf  \scriptfont\bffam=\sixbf
\scriptscriptfont\bffam=\fivebf\def\bf{\fam\bffam\ninebf}%
\tt\ttglue=.5em plus.25em minus.15em
\normalbaselineskip=11pt
\setbox\strutbox=\hbox{\vrule height8pt depth3pt width0pt}%
\let\sc=\ninesc\let\mc=\eightrm
\def\cyr{\ninecyr\cyracc}\def\cyri{\ninecyri\cyracc}
\let\big=\ninebig\normalbaselines\rm}

\def\eightpoint{\def\rm{\fam0\eightrm}%
  \textfont0=\eightrm \scriptfont0=\sixrm \scriptscriptfont0=\fiverm
  \textfont1=\eighti  \scriptfont1=\sixi  \scriptscriptfont1=\fivei
  \textfont2=\eightsy \scriptfont2=\sixsy \scriptscriptfont2=\fivesy
  \textfont3=\tenex   \scriptfont3=\tenex \scriptscriptfont3=\tenex
  \textfont\itfam=\eightit  \def\it{\fam\itfam\eightit}%
  \textfont\slfam=\eightsl  \def\sl{\fam\slfam\eightsl}%
  \textfont\ttfam=\eighttt  \def\tt{\fam\ttfam\eighttt}%
  \textfont\bffam=\eightbf  \scriptfont\bffam=\sixbf
  \normalbaselineskip=9pt
  \let\sc=\eightsc \let\mc=\sevenrm  
  \def\cyr{\eightcyr\cyracc}\def\cyri{\eightcyri\cyracc}
  \let\big=\eightbig  \normalbaselines\rm}%
\def\nospace{\nulldelimiterspace0pt\mathsurround0pt}%
\def\tenbig#1{{\hbox{$\left#1\vbox to8.5pt{}\right.\nospace$}}}%
\def\ninebig#1{{\hbox{$\textfont0=\tenrm\textfont2=\tensy
  \left#1\vbox to7.25pt{}\right.\nospace$}}}%
\def\eightbig#1{{\hbox{$\textfont0=\ninerm\textfont2=\ninesy
  \left#1\vbox to6.5pt{}\right.\nospace$}}}%

\def\nonextendedbold{
  \font\fiveb=cmb10 at 5pt
  \font\sixb=cmb10 at 6pt
  \font\sevenb=cmb10 at 7pt
  \font\eightb=cmb10 at 8pt
  \font\nineb=cmb10 at 9pt
  \font\tenb=cmb10
  \font\twelveb=cmb10 at 12pt
  \let\fivebf=\fiveb
  \let\sixbf=\sixb
  \let\sevenbf=\sevenb
  \let\eightbf=\eightb
  \let\ninebf=\nineb
  \let\tenbf=\tenb
  \let\twelvebf=\twelveb
}

\def\leftrighttop#1#2{
  \headline{\ifnum\pageno=1\hfil\else{\ninept #1 \hfil #2}\fi}
}

\def\firstnopagenum{
  \footline{\ifnum\pageno=1 \hfil \else \hfil{\rm \number\pageno}\hfil\fi}
}

\def\maketitle#1#2#3#4{
  \centerline {\titlefont #1}
  \medskip
  \centerline {\eightpt #2}
  \medskip
  \centerline {\tensc #3}
  \medskip
  \centerline {\tensc #4}
  \bigskip
}


\outer\def\floattext#1 #2. #3\par{
  $$
  \vbox{
    \hsize #1 true in
    \noindent{\bf #2.}\enskip #3
  }
  $$
}


\def\lsection#1\par{
  \bigskip\vskip\parskip
  \leftline{\sectionfont#1}\nobreak\medskip\noindent
}

\def\csection#1\par{
  \bigskip\vskip\parskip
  \centerline{\sectionfont#1}\nobreak\medskip\noindent
}

\def\rsection#1\par{
  \bigskip\vskip\parskip
  \rightline{\sectionfont#1}\nobreak\medskip\noindent
}
\def\section{\lsection}

\def\boldlabel#1. {\noindent{\bf #1.\enspace}}
\def\subsection#1. {\medskip\noindent{\bf #1.\enspace}}



\font\tenfrak=eufm10
\font\sevenfrak=eufm7
\font\fivefrak=eufm5
\newfam\frakfam
\textfont\frakfam=\tenfrak
\scriptfont\frakfam=\sevenfrak
\scriptscriptfont\frakfam=\fivefrak

\def\janksc#1#2 {#1{\eightpt#2}}
\def\jankscsp#1#2 {#1{\eightpt#2}\ }
\def\scproclaim#1.#2\par{\noindent\jankscsp #1.\enspace{\it#2\par}}

\def\bref#1{[#1]}
\def\ref#1{[#1]}

\def\quote{
  \begingroup
    \baselineskip 10pt
    \parfillskip 0pt
    \interlinepenalty 10000
    \leftskip 0pt plus 40pc minus \parindent
    \let\rm=\quoterm\let\sl=\quotesl\everypar{\sl}
    \obeylines
}
\def\author#1(#2){\nobreak\smallskip\rm--- \rm#1\unskip\enspace(#2)\par\endgroup}

\def\titlefont{\twelvebf}
\def\sectionfont{\tenssbx}
\def\quoterm{\eightssq}
\def\quotesl{\eightssqi}



\def\bookheader#1#2{
  \nopagenumbers
  \def\leftheadline{{\rm\folio}\hfil{\eightpoint#1}\hfil}
  \def\rightheadline{\hfil{\eightpoint#2}\hfil{\rm\folio}}
  \headline{\ifodd\pageno{\ifnum\pageno<2\hfil\else\rightheadline\fi}\else\leftheadline\fi}
}

\tenpoint



\def\xskip{\hskip 7pt plus 3pt minus 4pt}

\def\proof{\medbreak\noindent{\it Proof.}\xskip\ignorespaces}

\def\slug{\quad\hbox{\kern1.5pt\vrule width2.5pt height6pt depth1.5pt\kern1.5pt}\medskip}
\def\noskipslug{\quad\hbox{\kern1.5pt\vrule width2.5pt height6pt depth1.5pt\kern1.5pt}}

\newdimen\algindent
\newif\ifitempar \itempartrue 
\def\algindentset#1{\setbox0\hbox{{\bf #1.\kern.25em}}\algindent=\wd0\relax}
\def\algbegin #1 #2{\algindentset{#21}\alg #1 #2} 
\def\aalgbegin #1 #2{\algindentset{#211}\alg #1 #2} 
\def\alg#1(#2). {\medbreak 
  \noindent{\bf#1}({\it#2\/}).\xskip\ignorespaces}
\def\algstep#1.{\ifitempar\smallskip\noindent\else\itempartrue
  \hskip-\parindent\fi
  \hbox to\algindent{\bf\hfil #1.\kern.25em}%
  \hangindent=\algindent\hangafter=1\ignorespaces}



\def\ZZ{{\bf Z}}

\def\RR{{\bf R}}


\def\op#1{\mathop{\hbox{#1}}\nolimits}
\def\limitop#1{\mathop{\hbox{#1}}}


\def\pr{\op{\bf P}}



\def\eps{\epsilon}

\newcount\thmcount  
\thmcount=1
\newcount\sectcount  
\sectcount=1
\newcount\figcount  
\figcount=1
\newcount\eqcount  
\eqcount=1

\def\oldno#1{\eqno({\oldstyle#1})}
\def\refeq#1{({\oldstyle#1})}
\def\adveq{\oldno{\the\eqcount}\global\advance\eqcount by 1}  
\def\advthm{\the\thmcount\global\advance \thmcount by 1}

\def\advsect{\section\the\sectcount\global\advance\sectcount by 1. }

\def\caption#1{\centerline{\ninepoint{\bf Fig.~\the\figcount\global\advance\figcount by 1.\enspace}#1}}

\outer\def\parenproclaim #1 (#2).#3\par{\medbreak
  \noindent{\bf #1}\enspace\rm({\it #2\/}).\nobreak\ignorespaces{\sl #3\par}
  \ifdim\lastskip<\medskipamount \removelastskip\penalty55\medskip\fi}


\newdimen\axiomindent
\def\axiomindentset#1{\setbox0\hbox{{\bf #1.\kern.25em}}\axiomindent=\wd0\relax}
\def\axiom#1. [#2.]{\ifitempar\par\noindent\else\itempartrue
  \hskip-\parindent\fi%
  \hbox to\axiomindent{\bf\hfil #1.\kern.25em}%
  \hangindent=\axiomindent\hangafter=1[{\it #2.}]}

\input eplain

\def\bar{\overline}

\def\Eta{\op{\bf H}}
\def\Etaseven{\op{\sevenbf H}}

\def\AA{\op{\bf A}}
\def\MM{\op{\bf M}}
\def\ss{\op{\bf s}}
\def\FF{\op{\bf F}}

\def\given{\mathbin{|}}
\def\argmax{\limitop{\rm arg$\,$max}}

\def\bref#1{\special{ps:[/pdfm { /big_fat_array exch def big_fat_array 1 get 0
0 put big_fat_array 1 get 1 0 put big_fat_array 1 get 2 0 put big_fat_array pdfmnew } def}%
[\hlstart{name}{}{bib#1}#1\hlend]}

\magnification=\magstephalf
\hoffset=40pt \voffset=28pt
\hsize=29pc  \vsize=45pc  \maxdepth=2.2pt  \parindent=19pt
\nopagenumbers
\def\leftheadline{{\rm\folio}\hfil{\eightpoint MARCEL K. GOH}\hfil}
\def\rightheadline{\hfil{\eightpoint ON AN ENTROPIC ANALOGUE OF ADDITIVE ENERGY}\hfil{\rm\folio}}
\headline={\ifodd\pageno{\ifnum\pageno<2\hfil\else\rightheadline\fi}\else\leftheadline\fi}

\enablehyperlinks

\ifpdf
\hlopts{bwidth=0}
\pdfoutline goto name {intro} {Introduction}%
\pdfoutline goto name {energy} {Entropic additive energy}%
\pdfoutline goto name {large} {The large-energy regime}%
\pdfoutline goto name {small} {The small-energy regime}%
\pdfoutline goto name {mult} {Multiplicative energy and the sum-product problem}%
\pdfoutline goto name {acks} {Acknowledgements}%
\pdfoutline goto name {refs} {References}%
\else
\special{ps:[/PageMode /UseOutlines /DOCVIEW pdfmark}%
\special{ps:[/Count -0 /Dest (intro) cvn /Title (Introduction) /OUT pdfmark}%
\special{ps:[/Count -0 /Dest (energy) cvn /Title (Entropic additive energy) /OUT pdfmark}%
\special{ps:[/Count -0 /Dest (large) cvn /Title (The large-energy regime) /OUT pdfmark}%
\special{ps:[/Count -0 /Dest (small) cvn /Title (The small-energy regime) /OUT pdfmark}%
\special{ps:[/Count -0 /Dest (mult) cvn /Title (Multiplicative energy and the sum-product problem) /OUT pdfmark}%
\special{ps:[/Count -0 /Dest (acks) cvn /Title (Acknowledgements) /OUT pdfmark}%
\special{ps:[/Count -0 /Dest (refs) cvn /Title (References) /OUT pdfmark}%
\fi

\maketitle{On an entropic analogue of additive energy}{}{Marcel K. Goh}{\sl Department of
Mathematics and Statistics, McGill University}

\floattext4 \ninebf Abstract.
\ninepoint
Recent advances have linked various statements involving sumsets and cardinalities
with corresponding statements involving sums of random variables and entropies.
In this vein, this paper shows that the quantity
$2\Eta\{X, Y\} - \Eta\{X+Y\}$
is a natural entropic
analogue of the additive energy $E(A,B)$ between two sets. We develop some basic theory surrounding
this quantity, and demonstrate its role in the proof of Tao's entropy variant of the
Balog--Szemer\'edi--Gowers theorem. We examine the regime where entropic additive energy
is small, and discuss a family of random variables related to Sidon sets. In finite fields,
one can define an entropic multiplicative energy as well, and we formulate sum-product-type conjectures
relating these two entropic energies.
\smallskip
\noindent\boldlabel Keywords. Sumsets, entropy, additive energy, Sidon sets.
\smallskip
\noindent\boldlabel Mathematics Subject Classification. 11B13, 94A17.

\advsect Introduction
\hldest{xyz}{}{intro}

{\sc In a {\oldstyle 1960} article}
on metaphysics and mathematics~\bref{19}, the eminent number theorist and
algebraic geometer Andr\'e Weil wrote, ``{\it Rien n'est plus f\'econd,
tous les mathematiciens le savent, que ces obscures analogies, ces troubles reflets d'une th\'eorie
\`a une autre, ces furtives caresses, ces brouilleries inexplicables; rien aussi ne donne
plus de plaisir au chercheur.}\footnote{$^1$}{\eightpoint ``Nothing
is more fertile, all mathematicians know it, than these obscure analogies, these turbid reflections
from one theory to another, these furtive caresses, these inexplicable misunderstandings; nothing gives
more pleasure to the researcher.''}''
Over the last two decades or so, a small body of literature has emerged concerning
an obscure analogy between the worlds of additive combinatorics and information theory.

To go from the combinatorial setting to the information-theoretic one in this analogy
is as simple as performing the
following symbolic manoeuvre. One replaces sets $A$
with random variables $X$, (logarithms of) cardinalities $|A|$ with entropies $\Eta\{X\}$,
and sumsets $A+B$ with sums $X+Y$ of random variables.
The resulting scenario falls into one of three broad categories, as laid out in 2009 by
I.~Z.~Ruzsa~\bref{16}.
\medskip
\item{i)} The cardinality statement is provably equivalent to its entropy analogue.
\smallskip
\item{ii)} It is unknown if the statements are equivalent,
but the proof of one statement can be adapted to prove the other.
\smallskip
\item{iii)} The cardinality statement is provably true and the other remains a conjecture.
\medskip
Outside these three categories are examples of a cardinality statement being untrue despite
its entropy version being true
(see, e.g., the remark below Proposition~1.4 of~\bref{9})
and vice versa (see, e.g., Section~1.12 of~\bref{18}).
As such, the utility of the cardinality-to-entropy translation procedure might have seemed hopelessly tenuous
before the recent proof of the polynomial Freiman--Ruzsa conjecture
due to W.~T.~Gowers, B.~Green, F.~Manners, and T.~Tao (in characteristic $2$ first~\bref{6}, then in odd
characteristic~\bref{7}).

In this paper, we develop theory surrounding an entropic analogue of additive energy. The construction
appears in~\bref{6}, but no mention is made of its relation to additive energy. We justify the definition
by showing that it enjoys many properties analogous to those of combinatorial additive energy, and use these
properties to simplify some analysis of T.~Tao~\bref{18}
regarding the regime where entropic doubling is small (and thus entropic additive energy is large).
Then we perform a similar study of the situation in which the additive energy is small, and discuss
extremal cases related to Sidon sets. To close the paper, we consider an entropic multiplicative energy
as well, and formulate analogues of sum-product conjectures.

Very recently, a preprint of R.~Li, L.~Gavalakis, and I.~Kontoyiannis~\bref{13}
has generalised many of the results
in the present paper to differential entropy, obtaining improved bounds in certain instances.
Given that their paper references many of our propositions directly, in revising our paper we have chosen
to keep most of its content
unchanged to preserve the narrative thread of the literature.
However, at appropriate times
we shall make note of the advances made in the newer paper.

In the following section, we give a thorough and
unhurried introduction to Shannon entropy and its basic properties.
The reader who is already familiar with the rudiments of this theory may choose to skip to Section~3,
where the paper begins in earnest.

\advsect The Khintchine--Shannon axioms

This presentation of Shannon entropy and its properties draws
heavily from two lectures given by W.~T.~Gowers in~2020. Where a
proof in this section has been credited to someone else, this is the attribution supplied by Gowers.

All random variables in this paper are understood to be discrete.
The {\it entropy} $\Eta\{X\}$ of a random variable $X$ is a value in $[0,\infty]$ that represents
the amount of information that can be encoded in the random variable. For example, if $X$ is a constant
random variable, then $\Eta\{X\}$ should be zero, since one gains no information by learning the value
of $X$, and if $X$ is uniformly distributed on $\{0,1\}^n$, then we should like $X$ to be proportional
to $n$. The following axioms, which are known as the {\it Khintchine--Shannon axioms}, give rise to
a functional $\Eta$ that satisfies our desired properties.
\medskip
\item{a)} ({\it Invariance.}) If $X$ takes values in $A$, $Y$ takes values in $B$,
$\phi:A\to B$ is a bijection, and $\pr\{Y = \phi(a)\} = \pr\{X = a\}$ for all $a\in A$,
then $\Eta\{X\} = \Eta\{Y\}$.
\smallskip
\item{b)} ({\it Extensibility.}) If $X$ takes values in $A$ and $Y$ takes values in $B$ for
a set $B$ such that $A\subseteq B$, and furthermore $\pr\{Y=a\} = \pr\{X=a\}$ for all $a\in A$,
then $\Eta\{X\} = \Eta\{Y\}$.
\smallskip
\item{c)} ({\it Continuity.}) The quantity $\Eta\{X\}$ depends continuously on the probabilities
$\pr\{X=a\}$.
\smallskip
\item{d)} ({\it Maximisation.}) Over all possible random variables $X$ taking values in a finite
set $A$, the quantity $\Eta\{X\}$ is maximised for the uniform distribution.
\smallskip
\item{e)} ({\it Additivity.})
For $X$ taking values in $A$ and $Y$ taking values in $B$, we have the formula
$$ \Eta\{X,Y\} = \Eta\{ X \given Y\} + \Eta\{Y\},\adveq$$
where $\Eta\{X,Y\} = \Eta\bigl\{ (X,Y)\bigr\}$ and
$$ \Eta\{X\given Y\} = \sum_{y\in B} \pr\{Y=y\} \Eta\{ X\given Y = y\}.\adveq$$
\medskip
We shall take it on faith that there actually exists a functional on random variables satisfying these
axioms. In fact, there is only one such functional up to multiplication by a real
constant, so we shall add the following normalising axiom.
\medskip
\item{f)} ({\it Normalisation.}) If $X$ is uniformly distributed on $\{0,1\}$, then
$\Eta\{X\} = 1$.
\medskip
Notationally, one would expect that $\Eta\{X\given Y\} = \Eta\{X\}$ if $X$ and $Y$ are independent.
This is the first proposition we will carefully prove, using only the axioms.

\proclaim Proposition A. Let $X$ and $Y$ be independent random variables. Then
$\Eta\{X\given Y\} = \Eta\{X\}$ and consequently
$\Eta\{X,Y\} = \Eta\{X\} + \Eta\{Y\}$.

\proof Suppose $Y$ takes values in a set $B$. Then for all $y\in B$, the distribution of
$X$ and $X$ given that $Y = y$ is the same, so
$$\Eta\{X \given Y\} = \sum_{y\in B} \pr\{Y = y\} \Eta\{X\given Y = y\} = \sum_{y\in B}
\pr\{Y = y\} \Eta\{X\} =\Eta\{X\}.\adveq$$
The second version of the statement follows from the additivity axiom.\slug

We will sometimes use the notation $X^n$ to denote the vector $(X_1,\ldots,X_n)$ where the $X_i$
are independent copies of the random variable $X$. The following corollaries
are each proved by induction; the third of these is often known as the {\it chain rule}.

\edef\corxn{B}
\proclaim Corollary B. We have $\Eta\{X^n\} = n\Eta\{X\}$.\slug

\edef\cortwon{C}
\proclaim Corollary C. If $X$ is uniformly distributed on a set of size $2^n$, then
$$\Eta\{X\} = n.\noskipslug\adveq$$

\parenproclaim Corollary~{D} (Chain rule).
Let $X_1,\ldots,X_n$ be random variables. Then
$$\Eta\{X_1,\ldots,X_n\} = \Eta\{X_1\} + \Eta\{X_2\given X_1\}
+ \cdots+\Eta\{X_n \given X_1,\ldots,X_{n-1}\}.\noskipslug\adveq$$

Next we establish the intuitive statement
that the entropy of a uniform random variable supported on a set $A$ is at most
the entropy of a uniform random variable supported on a superset $B$ of $A$.

\edef\propunifineq{E}
\proclaim Proposition E. Let $A\subseteq B$ with $B$ finite,
let $X$ be uniformly distributed on $A$, and let $Y$ be uniformly distributed on $B$.
Then $\Eta\{X\} \le \Eta\{Y\}$, with equality if and only if $A = B$.

\proof By the extensibility axiom, $\Eta\{X\}$ is not affected if we regard $X$ as a function
taking values in $B$, and by the maximisation axiom, $\Eta\{X\}\le \Eta\{Y\}$, since $Y$ is
uniform on $B$.

If $A = B$, then it is clear that $\Eta\{X\} = \Eta\{Y\}$, since $X$ and $Y$
are the same random variable.

On the other hand, say $|A| = m$ and $|B| = n$ with $m<n$. We want to conclude the
strict inequality $\Eta\{X\} < \Eta\{Y\}$.
Pick $k$ such that $m^k \le n^{k-1}$, so that $A^k$ can be embedded (as a set) into $B^{k-1}$.
By invariance we can view $X^k$ as a random variable that is uniform on a subset of $B^{k-1}$.
so by the inequality we showed in the first paragraph of this proof, as well as Corollary~{\corxn},
we have
$$n\Eta\{X\} = \Eta\{X^n\} \le \Eta\{Y^{n-1}\} = (n-1)\Eta\{Y\},\adveq$$
whence $\Eta\{X\} < \Eta\{Y\}$.\slug

The argument we have just seen is an instance of the ``tensor power trick,''
in which one amplifies a known inequality by applying it to suitably high powers of the objects
involved.

\medskip\boldlabel Entropy and cardinality.
We often relate cardinalities of finite sets to entropies of random variables by invoking the
following proposition. Here and elsewhere in the paper we write $\log$ to denote the binary logarithm.
This base can be changed by modifying the normalisation axiom.

\edef\proplogeq{F}
\proclaim Proposition F. Let $X$ be a uniform random variable on a finite set $A$.
Then
$$\Eta\{X\} = \log |A|.\adveq$$

\proof For any positive integer $n$ we can let
$X^n$ denote a tuple of independent copies of $X$; Corollary~{\corxn} tells us that $\Eta\{X^n\} = n\Eta\{X\}$.
Let $m$ be such that $2^m \le |A|^n \le 2^{m+1}$ and hence
$${m\over n} \le \log |A| \le {(m+1)\over n}.\adveq$$
Let $Y$ be uniform on a set of size $2^m$,
and let $Z$ be uniform on a set of size $2^{m+1}$. Then by Corollary~{\cortwon} we have
$\Eta\{Y\} = m $ and $\Eta\{Z\} = (m+1)$. Then by Proposition~{\propunifineq}
we have
$${m\over n} \le \Eta\{X\} \le {(m+1)\over n}.\adveq$$
In other words, $\Eta\{X\}$ satisfies the same bounds as $\log |A|$. Taking $n$ large, we can make
these bounds arbitrarily tight, proving the claim.\slug

The maximisation axiom then gives the following corollary.

\proclaim Corollary G. Let $X$ be a random variable supported on a finite set $A$.
Then
$$\Eta\{X\} \le \log |A|.\noskipslug\adveq$$

\medskip\boldlabel Functions of random variables.
If $Y$ is a random variable such that $Y = f(X)$ for some random variable $Y$ and some function $f$,
then we say that $Y$ is {\it determined by} $X$ or $X$ {\it determines} $Y$.
We want to show that $\Eta\{Y\}\le \Eta\{X\}$, which reflects the idea that we get more information
from $X$ than from $Y$. This seems to require a couple of steps, though the lemmas we prove along the way
are useful in their own right.

\edef\lemdeterminesbijection{H}
\proclaim Lemma H. If $Y = f(X)$ then $\Eta\{X\} = \Eta\{Y\} + \Eta\{X\given Y\}$.

\proof There is a a bijection between values $x$ taken
by $X$ and values $\bigl(x,f(x)\bigr)$ taken by $(X,Y)$, so we have
$$\Eta\{X\} = \Eta\{X,Y\} = \Eta\{Y\} + \Eta\{X\given Y\}\adveq$$
by invariance and additivity.\slug

We are now done if we can show that entropy is nonnegative.
This is a corollary of the following lemma, whose proof is a modification of
one due to S.~Eberhard.

\proclaim Proposition I. Let $X$ be a discrete random variable supported on a set $A$ and let
$$a^* = \argmax_{a\in A} \pr\{X = a\}.$$
Then
$$\pr\{X = a^*\} \ge 2^{-\Etaseven\{X\}}.\adveq$$

\proof First we will work in the case where there exists $n$ such that $\pr\{X = a\}$ is a multiple of
$1/n$ for all $a\in A$. Let $Y$ be uniformly distributed on $[n]$ and let $\{E_a\}_{a\in A}$
be a partition of $[n]$ such that $|E_a| = n\pr\{X= a\}$ for all $a\in A$, and let $Z = a$
if $Y\in E_a$. This definition makes $Z$ and $X$ identically distributed, so $\Eta\{Z\} = \Eta\{X\}$
by the invariance axiom.

For every $a\in A$, the conditional entropy $\Eta\{Y \given Z = a\}$ is uniformly
distributed on a set of size $|E_a|$. From our choice of $a^*$
we have $|E_{a^*}| \ge |E_a|$ for all $a\in A$. Hence by Proposition~{\proplogeq}, we have
$$\Eta\{Y\given Z\} = \sum_{a\in A} \pr\{X = a\} \Eta\{Y\given X = a\}
= \sum_{a\in A} \pr\{X = a\} \log |E_a| \le \log |E_{a^*}|.\adveq$$
Since $Z$ is determined by $Y$, we have $\Eta\{Y\} = \Eta\{Z\} + \Eta\{Y\given Z\}$ by the previous lemma,
and by another invocation of Proposition~{\proplogeq}, we have
$$\eqalign{
\Eta\{Z\} &= \Eta\{Y\} - \Eta\{Y\given Z\} \cr
&\ge \log n - \log |E_{a^*}| \cr
&\ge \log \biggl({n\over |E_{a^*}|}\biggr) \cr
&= \log\biggl({1\over \pr\{X = a^*\}}\biggr), \cr
}\adveq$$
and hence $2^{-\Etaseven\{X\}} = 2^{-\Etaseven\{Z\}} \le \pr\{X = a^*\}$.

The general case follows from the continuity axiom.\slug

In this proof, we came dangerously close to obtaining the formula for entropy.
To emphasise the axiomatic approach,
we shall refrain from writing this formula out explicitly, but a keen reader might enjoy working
through the details of its derivation.

From the fact that $\pr\{X = a^*\} \le 1$, we immediately conclude that entropy is nonnegative.

\proclaim Corollary J. Let $X$ be a discrete random variable taking values in a finite
set $A$. Then $\Eta\{X\} \ge 0$.\slug

With this observation (and Lemma~{\lemdeterminesbijection}), we see that a
random variable has a smaller entropy than one by which it is determined.

\edef\cordetermines{K}
\proclaim Corollary K. If $Y = f(X)$ then $\Eta\{X\} \ge \Eta\{Y\}$.\slug

Next we show that a random variable has zero entropy if and only if it is constant.
This reflects the idea that the variables from which we get no information are those which take
the same value no matter what.

\proclaim Proposition L. Let $X$ be a discrete random variable.
Then $\Eta\{X\} = 0$ if and only if it takes exactly one value.

\proof First suppose that $X$ takes only one value.
Let $a$ be the value of $X$ such that $\pr\{X=a\} = 1$.
Then $(X,X)$ equals $(a,a)$ with probability $1$ as well,
so $\Eta\{X\} = \Eta\{X,X\}$ by the invariance axiom. But it can easily be checked that
$X$ and $(X,X)$ are independent (we have
$$\pr\bigl\{X = a, (X,X)= (a,a)\bigr\} = \pr\{X = a\}\pr\bigl\{(X,X) = (a,a)\bigr\}\adveq$$
for instance), so $\Eta\{X,X\} = 2\Eta\{X\}$. Thus we conclude
that $\Eta\{X\} = 0$.

Now suppose that $X$ takes more than one value; let $A$ be the set of $a$ such that
$\pr\{X = a\} > 0$ and let $\alpha = \max_{a\in A} \pr\{X = a\}$. For all $n$ let $X^n$ denote
the tuple of $n$ independent copies of $X$; the maximum probability of any particular value (in $A^n$)
that $X^n$ takes is $\alpha^n$. But $\alpha < 1$ since $X$ takes more than one value, so for any
$\eps > 0$ we can find $n$ such that $\alpha^n < \eps$. This means that we can partition $A^n$ into
two disjoint sets $E$ and $F$ such that $\pr\{X^n \in E\}$ and $\pr\{X^n \in F\}$ are both
in the range $[1/2-\eps, 1/2+\eps]$.

Let $Y$ be the random variable taking the value $0$ if $X^n\in E$ and $1$ if $X^n\in F$.
By Corollary~{\corxn}, $\Eta\{X^n\} = n\Eta\{X\}$, and since $X^n$ determines $Y$,
$$\Eta\{X^n\} = \Eta\{Y\} + \Eta\{X^n\given Y\} \ge \Eta\{Y\}.\adveq$$
But $\Eta\{Y\}$ can be taken as close to $1$ as we like by choosing $\eps$ small enough,
by the normalisation and continuity axioms. So
$\Eta\{X\} \ge \Eta\{Y\}/n > 0$.\slug

\medskip\boldlabel Subadditivity and submodularity.
We now prove the intuitive fact that conditioning
cannot increase the entropy of a random variable.
The proof we present is due to C.~West.

\proclaim Proposition M. Let $X$ and $Y$ be discrete random variables. Then
$$\Eta\{X\given Y\} \le \Eta\{X\}.\adveq$$

\proof Let $A$ be the support of $X$ and $B$ be the support of $Y$.
First we consider the case that $X$ is uniform on $A$ (so $A$ is finite). Then by the definition
of conditional entropy,
$$\Eta\{X\given Y\} = \sum_{b\in B} \pr\{Y = b\} \Eta\{X \given Y= b\}.\adveq$$
But for each $b$, the random variable $(X\given Y = b)$ takes values in $A$, so its entropy is
bounded above by $\Eta\{X\}$ by the maximisation axiom. Hence the claim is true in this special case.

Next, suppose that $A$ and $B$ are both finite and suppose further that $\pr\{Y = b\}$
is rational for all $b$. Then there is an integer $n$ and integers $\{m_b\}_{b\in B}$
such that $\pr\{Y= b\} = m_b/n$ for all $b\in B$. Now partition $[n]$ into sets $\{E_b\}_{b\in B}$,
where $|E_b| = m_b$ for all $b\in B$. We define a random variable $Z$ by sampling uniformly at
random from $E_b$ if $Y = b$, and doing so independently of $(X\given Y=b)$. The result is
a random variable $Z$ that is uniform on $[n]$ and which is independent of $X\given Y$.
Furthermore, since $Z$ determines $Y$, we have $\Eta\{Z\} = \Eta\{Y,Z\}$ by the invariance
axiom. Hence
$$\eqalign{
\Eta\{X\given Y\} &= \Eta\{X\given Y, Z\} \cr
&= \Eta\{X,Y,Z\} - \Eta\{Y,Z\} \cr
&= \Eta\{X,Z\} - \Eta\{Z\} \cr
&= \Eta\{X\given Z\} \cr
&\le \Eta\{X\}, \cr
}\adveq$$
where the inequality on the last line follows from the previous paragraph.

The general case follows from the continuity axiom and the fact that any discrete random
variable, regarded as a vector in $l_1(\RR)$,
by vectors with finite support, and these in turn can be approximated by vectors of finite support
and rational coordinates.\slug

By the additivity axiom, the previous proposition is equivalent to the following subadditive property
of entropy:
$$\Eta\{X,Y\} \le \Eta\{X\} + \Eta\{Y\}\adveq$$
We shall use this to prove the important submodularity inequality.

\parenproclaim Proposition N (Submodularity).
Suppose $X$, $Y$, $Z$, and $W$ are random variables such that $(Z,W)$ determines
$X$, $Z$ determines $Y$, and $W$ also determines $Y$. Then
$$\Eta\{X\} + \Eta\{Y\} \le \Eta\{Z\} + \Eta\{W\}.\adveq$$

\proof Since $(Z,W)$ determines $X$, for each possible value $y$ that $Y$ takes,
$(Z,W\given Y= y)$ determines $(X\given Y=y)$. Summing over the support of $Y$ and
applying Corollary~{\cordetermines} gives us
$$ \Eta\{X\given Y\} \le \Eta\{Z,W \given Y\} \le \Eta\{Z\given Y\} + \Eta\{W\given Y\}.$$
Expanding the conditional entropies and rearranging yields
$$\Eta\{X,Y\} + \Eta\{Y\} \le \Eta\{Z,Y\} + \Eta\{W,Y\}.$$
But $Z$ and $W$ each determine $Y$, so $\Eta\{Z,Y\} = \Eta\{Z\}$ and
$\Eta\{W,Y\} = \Eta\{W\}$. Hence
$$\Eta\{X\} + \Eta\{Y\} \le \Eta\{X,Y\} + \Eta\{Y\} \le \Eta\{Z\} + \Eta\{W\},\adveq$$
as desired.\slug

\medskip\boldlabel Group-valued random variables.
If two random variables under consideration take values in the same abelian group,
as will often be the case for us, it makes sense to speak of entropies of sums
or differences. Since conditioning does not increase entropy, we have
$$\Eta\{X+Y\} \ge \Eta\{X+Y\given X\} = \Eta\{Y\given X\},\adveq$$
and the same thing with $X+Y$ replaced with $X-Y$ or $X$ replaced with $Y$. Consequently, if
$X$ and $Y$ are independent, then
\edef\eqindepmax{\the\eqcount}
$$\min\bigl(\Eta\{X+Y\}, \Eta\{X-Y\}\bigr) \ge \max\bigl(\Eta\{X\}, \Eta\{Y\}\bigr).\adveq$$

\medskip\boldlabel Conditionally independent trials.
Let $X$, $Y$, and $Z$ be random variables (among which no specific independence relations are imposed).
We say that
$X_1$ and $Y_1$ are {\it conditionally independent trials of $X$ and $Y$ relative to $Z$} if for all
$z$ in the range
of $Z$, the random variables distributed as $(X_1\given Z=z)$ and $(Y_1\given Z=z)$ are independent,
$(X_1\given Z=z)$ has the same distribution as $(X\given Z=z)$,
and similarly for $Y_1$ and $Y$.
In particular, if $X=Y$ and $X_1$ and $X_2$ are conditionally independent trials of $X$ relative
to $Z$, we have
$$\Eta\{X_1,X_2\given Z\} = \Eta\{X_1\given Z\} + \Eta\{X_2\given Z\} = 2\Eta\{X\given Z\},$$
by additivity and independence. From this we obtain
\edef\eqcondindep{\the\eqcount}
$$\Eta\{X_1,X_2,Z\} = 2\Eta\{X\given Z\} + \Eta\{Z\} = 2\Eta\{X,Z\} - \Eta\{Z\}.\adveq$$
Observe also that $(X_1,Z)$ and $(X_2,Z)$ both have the same distributions
as $(X,Z)$.

\advsect Entropic additive energy
\hldest{xyz}{}{energy}

For subsets $A$ and $B$ of the same abelian group, it stands to reason
that $|A+B|$ is much smaller than $|A|\cdot|B|$ if and only if there is
significant redundancy in $A+B$; that is, many elements of $A+B$ can
be expressed as $a+b$ in different ways. To capture this notion, it is customary to speak of the
quantity
$$E(A,B) = \bigl| \bigl\{ (a,a',b,b') \in A\times A\times B\times B : a + b = a'+b'\bigr\}\bigr|$$
as the {\it additive energy} between the two sets $A$ and $B$.

We now give an entropic version of $E(A,B)$.
Let $X$ and $Y$ be discrete random variables taking values
in the same
abelian group. Let $(X_1, Y_1)$ and $(X_2, Y_2)$ be conditionally independent trials of $(X,Y)$ relative to
$X+Y$. Thus $X_1 + Y_1 = X_2 + Y_2 = X+Y$.
Define the {\it entropic additive energy} between $X$ and $Y$ to be
$$ \AA\{X,Y\} = \Eta\{X_1, Y_1, X_2, Y_2\}.$$
By conditional independence and the fact that $(X_1, Y_1, X_2, Y_2)$ determines $X_1+Y_1 = X+Y$,
we can rewrite
$$\eqalign{
\AA\{X,Y\} &= \Eta\{X_1, Y_1, X_2, Y_2, X+Y\} \cr
&= 2\Eta\{X,Y,X+Y\} - \Eta\{X+Y\} \cr
&= 2\Eta\{X,Y\} - \Eta\{X+Y\},\cr
}$$
where in the second equality we applied~\refeq{\eqcondindep}. This second formula for
$\AA\{X,Y\}$ is the one we use more often, since it sidesteps the construction involving
the variables $(X_1, Y_1)$ and $(X_2,Y_2)$.
The construction above was employed in~\bref{6} but the quantity $\Eta\{X_1, Y_1, X_2, Y_2\}$
was not given any specific name in that paper. We believe it to be the most natural
analogue of additive energy in the entropy setting, due to its construction as well as the following
quantitative observation.

\edef\propbasicbounds{\the\thmcount}
\proclaim Proposition {\advthm}. Let $A$ and $B$ be subsets of an abelian group, and let $U_A$ and
$U_B$ be independent indicator random variables on $A$ and $B$, respectively. Then
$$\log{|A|^2 |B|^2\over |A+B|} \le \AA\{U_A, U_B\} \le \log E(A,B).\adveq$$

\proof By independence and the fact that $U_A + U_B$ is supported on $A+B$, we have
$$\eqalign{
\AA\{U_A, U_B\} &= 2 \Eta\{U_A, U_B\} - \Eta\{U_A + U_B\} \cr
&\ge 2\log |A| + 2\log |B| - \log |A+B| \cr
&= \log {|A|^2 |B|^2 \over |A+B|}, \cr
}\adveq$$
establishing the first inequality.

For the second inequality, we use the conditional definition of $\AA\{X,Y\}$.
In the construction above, the tuple $(X_1, Y_1, X_2, Y_2)$ takes values in the set
$$\bigl\{ (a,a',b,b') \in A\times A\times B\times B : a + b = a'+b'\bigr\},$$
so $\AA\{X,Y\} \le \log E(A,B)$.\slug

As a byproduct we have furnished an information-theoretic proof of the elementary inequality
$${|A|^2|B|^2\over |A+B|} \le E(A,B),$$
which can also be proved using a counting argument
and the Cauchy--Schwarz inequality (see Lemmas~2.9 and~2.10
in~\bref{17}).
On the subject of the
Cauchy--Schwarz inequality, it turns out that the inequality
$$E(A,B) \le E(A,A)^{1/2} E(B,B)^{1/2}\adveq$$
for the classical additive energy does {\it not} remain true when transferred to the entropy setting.

First we establish the shorthand
\edef\eqenergyonesetdef{\the\eqcount}
$$\AA\{X\} = \AA\{X,X'\} = 2\Eta\{X,X'\} - \Eta\{X+X'\} = 4\Eta\{X\} - \Eta\{X+X'\},\adveq$$
where $X'$ is an independent copy of $X$.
The claim is that one can find random variables $X$ and $Y$ with
\edef\eqcscounterexample{\the\eqcount}
$$\AA\{X,Y\} > {1\over 2} \AA\{X\} + {1\over 2} \AA\{Y\}.\adveq$$
To this end, let
$$A = \{-7,-5,-4,-3,0,4,5,7\},$$
considered as subset of the integers. (This is the smallest subset of the integers satisfying
$|A+A| > |A-A|$, as shown by P.~Hegarty in~\bref{10}.) We shall set $X = U_A$
and $Y = U_{-A}$, such that $X$ and $Y$ are independent. We clearly have
$\Eta\{X\} = \Eta\{Y\} = 3$, and it can easily be verified computationally that
$$\Eta\{X+Y\} \approx 4.507 \qquad\hbox{and}\qquad
\Eta\{X+X'\} = \Eta\{Y + Y'\} \approx 4.513 ,\adveq$$
where $X'$ is an independent copy of $X$ and $Y'$ is an independent copy of $Y$.
Since $X$ and $Y$ are independent, $\AA\{X,Y\} = 12 - \Eta\{X+Y\} \approx 7.493$,
while on the other hand, $\AA\{X\} = \AA\{Y\} = 12 - \Eta\{X+X'\} \approx 7.487$.
Hence $X$ and $Y$ satisfy~\refeq{\eqcscounterexample}.

\advsect The large-energy regime
\hldest{xyz}{}{large}

The additive energy $E(A,B)$ is bounded from above by
$|A|\cdot|B|\min\bigl(|A|, |B|\bigr)$, since if $a+b = a'+b'$, then any three terms determine the fourth.
Therefore, an additive energy on the order of $|A|^{3/2} |B|^{3/2}$ may be considered large.
Similarly, we have
$$\eqalign{
\AA\{X,Y\} &= 2\Eta\{X,Y\} - \Eta\{X+Y\} \cr
&\le 2\Eta\{X,Y\} - {1\over 2} \Eta\{X+Y\given X\} - {1\over 2} \Eta\{X+Y\given Y\} \cr
&\le 2\Eta\{X,Y\} - \min\bigl( \Eta\{Y\given X\}, \Eta\{X\given Y\}\bigr), \cr
&\le \Eta\{X,Y\} + \min\bigl( \Eta\{X\}, \Eta\{Y\}\bigr), \cr
}\adveq$$
so we might say that $\AA\{X,Y\}$ is large if it differs
from $(3/2)(\Eta\{X\} + \Eta\{Y\})$ by some constant error term.

In this section we show that $\AA\{X,Y\}$ being large corresponds
to the case in which $\Eta\{X+Y\}$ is small compared to $\Eta\{X\}$ and $\Eta\{Y\}$.
To our knowledge, the fact that large $\AA\{X,Y\}$ implies small $\Eta\{X+Y\}$ does not appear
in the literature. A stronger statement, that large $\AA\{X,Y\}$ implies small $\Eta\{X'+Y'\}$,
where $X'$ and $Y'$ are {\it independent} conditionings of $X$ and $Y$, is given by Tao's entropic
Balog--Szemer\'edi--Gowers theorem~\bref{18}, of which we shall record a simpler version
written in terms of the entropic additive energy.

Our starting point is this next
proposition from the world of sumsets, which states that small sumset implies large additive energy.
It follows directly from the Cauchy--Schwarz inequality.

\edef\propinversebalog{\the\thmcount}
\proclaim Proposition \advthm.
Let $A$ and $B$ be finite subsets of an abelian group.
If $|A+B|\le C|A|^{1/2} |B|^{1/2}$ for some constant $C$, then we have
$$E(A,B) \ge {1\over C} |A|^{3/2} |B|^{3/2}.\noskipslug\adveq$$

Somewhat surprisingly, if we convert
this statement into its entropic analogue in a na\"\i ve manner,
the implication goes the other way.
We do have a weak equivalence (with worse constants in one direction)
under the further assumption that the random variables in
question are not too dependent.

\edef\propnaive{\the\thmcount}
\proclaim Proposition \advthm. Let $X$ and $Y$ be discrete random variables taking values in the
same abelian group. If
\global\edef\eqenergybound{\the\eqcount}
$$\AA\{X,Y\} \ge {3\over 2} \Eta\{X\} + {3\over 2}\Eta\{Y\} - \log C,\adveq$$
for some constant $C$, then
\global\edef\eqsumbound{\the\eqcount}
$$\Eta\{ X+Y\} \le {1\over 2} \Eta\{X\} + {1\over 2} \Eta\{Y\} + \log C.\adveq$$
If one adds the further assumption that
$\Eta\{X,Y\} \ge \Eta\{X\} + \Eta\{Y\} - C'$,
then~{\rm \refeq{\eqsumbound}} implies~{\rm \refeq{\eqenergybound}} with a worse constant, namely,
we may only conclude
\global\edef\eqworsebound{\the\eqcount}
$$\AA\{X,Y\} \ge {3\over 2} \Eta\{X\} + {3\over 2}\Eta\{Y\} - \log C - 2C'.\adveq$$
In particular, if $X$ and $Y$ are independent, then we can recover~{\rm \refeq{\eqenergybound}}
from~{\rm \refeq{\eqsumbound}}.

\proof Assuming the lower bound on the additive energy, we have
$$2\Eta\{X,Y\} - \Eta\{X+Y\} \ge {3\over 2} \Eta\{X\} + {3\over 2} \Eta\{Y\} - \log C,\adveq$$
so
$$\eqalign{
\Eta\{X+Y\} &\le 2\Eta\{X,Y\} - {3\over 2} \Eta\{X\} - {3\over 2} \Eta\{Y\} + \log C \cr
&\le 2\Eta\{X\} + 2\Eta\{Y\} - {3\over 2} \Eta\{X\} - {3\over 2} \Eta\{Y\} + \log C \cr
&= {1\over 2} \Eta\{X\} + {1\over 2} \Eta\{Y\} + \log C .\cr
}\adveq$$
On the other hand, assuming this upper bound on $\Eta\{X+Y\}$, we have
$$\eqalign{
\AA\{X,Y\} &= 2\Eta\{X,Y\} - \Eta\{X+Y\}  \cr
&\ge 2\Eta\{X,Y\} - {1\over 2} \Eta\{X\} - {1\over 2} \Eta\{Y\} - \log C ,\cr
}\adveq$$
and if $2\Eta\{X,Y\} \ge 2 \Eta\{X\} + 2\Eta\{Y\} - 2C'$, then~\refeq{\eqworsebound} follows
directly.\slug

At this point the reader may be baffled that the implication apparently goes the wrong way in the entropic
version of the statement.
Let us get to the bottom of this.
If $A$ and $B$ are subsets of a finite abelian group, we can let $X$ and $Y$ be the uniform
distributions on $A$ and $B$, respectively. We have been operating under the belief that $\Eta\{X+Y\}$
should correspond (up to taking powers or logarithms) to the size of $A+B$. But this is not true, since
$X$ and $Y$ may be given a joint distribution that is not uniform on $A\times B$, even if
its marginals are uniform on $A$ and $B$.

For example, let $A$ and $B$ be subsets of $G$ and consider any regular bipartite graph $H$ on the
vertex set $A\cup B$. Let $(X,Y)$ be defined by sampling an edge from $H$ uniformly at random, letting
$X$ be its endpoint in $A$ and $Y$ be its endpoint in $B$. Since the graph is regular, $X$ is uniform
on $A$ and $Y$ is uniform on $B$, but $X+Y$ can only take values $a+b$ where $(a,b)$ is an edge of $H$.
In other words, $X+Y$ samples from the {\it partial sumset}
$$ A +_H B = \bigl\{ a + b : (a,b)\in E(H)\bigr\},$$
where elements that are represented more times as the sum of edge endpoints are given a greater weight.
The way to properly recover the ordinary sumset $A+B$ is to let $H$ be all of $A\times B$, in which case
$X$ and $Y$ are independent. The extra assumption we added in Proposition~{\propnaive} is analogous to
stipulating that $|H| \ge C |A|\cdot |B|$, so that the resulting $X$ and $Y$ are ``nearly'' independent.

Simply put, in the entropic setting the bound~\refeq{\eqenergybound} is
stronger than~\refeq{\eqsumbound} because the entropy $\Eta\{X,Y\}$ of the joint distribution appears in
the formula for $\AA\{X,Y\}$ and thus forced to be somewhat large,
whereas~\refeq{\eqsumbound} says nothing whatsoever about this joint distribution and if
$X$ and $Y$ are rather dependent, its entropy could be small.

\medskip\boldlabel The Balog--Szemer\'edi--Gowers theorem.
The converse to Proposition~{\propinversebalog} does not hold in general; that is, large additive energy
does not necessarily imply a small sumset. However, there does exist a partial converse,
which says that if sets $A$ and $B$ have a large additive energy, then there are dense subsets $A'\subseteq A$
and $B'\subseteq B$ such that the sumset $|A+B|$ is small. This is the celebrated
Balog--Szemer\'edi--Gowers theorem, proven in 1994 by A.~Balog and E.~Szemer\'edi
in~\bref{1}, then subsequently improved upon in 1998 by W.~T.~Gowers~\bref{5}.

\parenproclaim Theorem {\advthm} (Balog--Szemer\'edi--Gowers theorem). Let $A$ be a finite
subset of an abelian group with $E(A,B) \ge c |A|^{3/2} |B|^{3/2}$. Then there are subsets
$A'\subseteq A$ and $B'\subseteq B$ with $|A'| \ge c' |A|$ and $|B'|\ge c''|B|$ such that
$$|A'+B'| \le C |A|^{1/2}|B|^{1/2},\adveq$$
where $c'$, $c''$, and $C$ depend only on $c$.\slug

We shall now present a leisurely exposition of Tao's entropic analogue of the Balog--Szemer\'edi--Gowers
theorem~\bref{18}, which was improved upon in~\bref{6}.
Our object is to reframe its statement and proof in the
the language of entropic additive energy, illustrating more directly its link to
the classical energy formulation of the Balog--Szemer\'edi--Gowers theorem above. No originality
is claimed for the proofs in this subsection.

In the entropy setting, the operation on random variables $X$ and $Y$ that corresponds to the unrestricted
sumset is the {\it independent sum} $\Eta\{X'+Y'\}$, rather than the ordinary sumset $\Eta\{X+Y\}$,
which may correspond to some incomplete bipartite graph, as we saw above.
The operation on random variables that corresponds to taking subsets is
conditioning. (As a sanity check, recall that conditioning never increases entropy, just
as taking subsets never increases cardinality.) The Balog--Szemer\'edi--Gowers theorem gives
us subsets between which we can take a {\it bona fide} sumset, so its entropic analogue
should return conditionings $X'$ and $Y'$ of $X$ and $Y$ relative to some random variable $Z$,
such that
\medskip
\item{i)} $X'$ and $Y'$ are conditionally independent relative to $Z$;
\smallskip
\item{ii)} the entropies $\Eta\{X'\given Z\}$ and
$\Eta\{Y'\given Z\}$ are not too small compared to their unconditioned analogues; and
\smallskip
\item{iii)} $\Eta\{X'+Y'\given Z\}$ is small.
\medskip
In fact, the conditioning we shall perform is exactly the one used to define additive energy.

First we state an inequality that appears within a proof in~\bref{6},
We have rewritten its statement in terms of the additive energy and excised the relevant portion
of the proof.

\edef\lematwo{\the\thmcount}
\parenproclaim Lemma {\advthm} ({\rm\bref{6},} Lemma A.2).
Let $X$ and $Y$ be discrete random variables taking values in the same abelian
group. Let $(X_1, Y_1)$ and $(X_2,Y_2)$ be conditionally
independent trials of $(X,Y)$ relative to $X+Y$. Then we have
$$\max\bigl(\Eta\{X_1- X_2\}, \Eta\{X_1-Y_2\}\bigr) \le 2\Eta\{X\} + 2\Eta\{Y\} - \AA\{X,Y\}.\adveq$$

\proof First we perform the proof for $X_1 - Y_2$.
Submodularity gives us
$$\Eta\{X_1, Y_1, X_1-Y_2\} + \Eta\{X_1-Y_2\} \le
\Eta\{ X_1, X_1-Y_2\} + \Eta\{Y_1, X_1-Y_2\}.\adveq$$
Since $X_1 + Y_1 = X+Y = X_2 + Y_2$, given $(X_1, Y_1, X_1-Y_2)$ we can recover
the values of $X_2$ and $Y_2$. So $(X_1, Y_1, X_1-Y_2)$ and $(X_1,Y_1,X_2,Y_2)$ determine each
other and hence
$$\Eta\{X_1, Y_1, X_1-Y_2\} = \Eta\{X_1, Y_1, X_2, Y_2\} = 2\Eta\{X,Y\} - \Eta\{X+Y\}.\adveq$$
On the other side of the inequality, we have
$$\Eta\{ X_1, X_1-Y_2\} = \Eta\{X_1, Y_2\} \le \Eta\{X\} + \Eta\{Y\},\adveq$$
and similarly
$$\Eta\{ Y_1, X_1-Y_2\} = \Eta\{Y_1, X_2-Y_1\} = \Eta\{X_2, Y_1\} \le \Eta\{X\} + \Eta\{Y\}.\adveq$$
Therefore,
$$\eqalign{
\Eta\{X_1-Y_2\} &\le \Eta\{X+Y\} - 2\Eta\{X,Y\} + 2\Eta\{X\} + 2\Eta\{Y\} \cr
&= 2\Eta\{X\} + 2\Eta\{Y\} - \AA\{X,Y\}.\cr
}\adveq$$
The same holds with $Y_2$ replaced by $X_2$.\slug

The bound on $\Eta\{X_1-Y_2\}$ is important to the proof of the polynomial Freiman--Ruzsa theorem over
$\FF_2^n$,
but it is the bound on $\Eta\{X_1-X_2\}$ that we require in the next proof.
As alluded to in~\bref{6}, this
lemma allows us to prove a better version of the entropic Balog--Szemer\'edi--Gowers
given as Theorem~3.1 of~\bref{18}. The proofs are slightly different because
a different conditioning is performed in that paper.
The two hypotheses in Tao's version of the statement
are
$$\Eta\{X,Y\} \ge \Eta\{X\} + \Eta\{Y\} - \log C$$
and
$$\Eta\{X+Y\} \le {1\over 2}\Eta\{X\} + {1\over 2}\Eta\{Y\} + \log C,$$
but by Proposition~{\propnaive} with $C'$ set to $\log C$, we may replace both of these
hypotheses with the single hypothesis
$$\AA\{X,Y\}\ge {3\over 2} \Eta\{X\} + {3\over 2}\Eta\{Y\} - 3\log C.$$
For cleanliness, however, by appropriately modifying $C$ we
shall simply dispense with the factor of $3$.

\parenproclaim Theorem {\advthm} (Entropic Balog--Szemer\'edi--Gowers theorem).
Let $X$ and $Y$ be discrete random variables taking values in the same abelian group, and
suppose that
$$\AA\{X,Y\} \ge {3\over 2} \Eta\{X\} + {3\over 2} \Eta\{Y\} - \log C$$
for some constant $C$.
Then letting $(X_1, Y_1)$ and $(X_2,Y_2)$ be conditionally
independent trials of $(X,Y)$ relative to $X+Y$, we have
$$ \Eta\{ X_1 \given X+Y\} \ge \Eta\{X\} - 2\log C\adveq$$
and
$$ \Eta\{ Y_2 \given X+Y\} \ge \Eta\{Y\} - 2\log C.\adveq$$
Furthermore, the variables $X_1$ and $Y_2$ are conditionally independent relative to $X+Y$, and
we have
$$\Eta\{ X_1 + Y_2 \given X+Y\} \le {1\over 2}\Eta\{X\} + {1\over 2} \Eta\{Y\} + \log C.\adveq$$

\proof Using the coupling $X+Y = X_1+Y_1 = X_2 + Y_2$, we have
$$\eqalign{
\Eta\{X_1\given X + Y\} &= \Eta\{X_1, X_1+Y_1\} - \Eta\{X+Y\} \cr
&= \Eta\{X,Y\} - \Eta\{X+Y\} \cr
&= \AA\{X,Y\} - \Eta\{X,Y\} \cr
&\ge {3\over 2} \Eta\{X\} + {3\over 2} \Eta\{Y\} - \log C - \Eta\{X,Y\} \cr
&\ge {1\over 2} \Eta\{X\} + {1\over 2}\Eta\{Y\} -\log C.\cr
}\adveq$$
where in the fourth line we used the hypothesis on $\AA\{X,Y\}$, and in the last line we observed
that $\Eta\{X\} + \Eta\{Y\} - \Eta\{X,Y\} \ge 0$. The bound
$$\Eta\{Y_2, X+Y\} \ge {1\over 2} \Eta\{X\} + {1\over 2}\Eta\{Y\} -\log C\adveq$$
is shown in the exact same way; only the first step differs.

Now, taking the sum of both these bounds, we arrive at
$$\Eta\{X_1\given X+Y\} + \Eta\{Y_2 \given X+Y\} \ge \Eta\{X\} + \Eta\{Y\} - 2\log C.\adveq$$
From this one deduces
$$\eqalign{
\Eta\{X_1\given X+Y\} &\ge \Eta\{X\} + \Eta\{Y_2\} - \Eta\{Y_2\given X+Y\}  - 2\log C \cr
&\ge \Eta\{X\} - 2\log C.\cr
}\adveq$$
The corresponding lower bound on $\Eta\{Y_2\given X+Y\}$ is proved similarly.

It remains to prove the upper bound on $\Eta\{X_1 + Y_2 \given X+Y\}$.
Note that $(X_1, Y_2, X+Y)$ and $(X_1 - X_2, X+Y)$ jointly determine $(X_1, X_2, X+Y)$.
Then given $X_1 - X_2$ and $X+Y$ we can calculate
$$X_1 + Y_2 = X_1 - X_2 + X_2 + Y_2 = (X_1 - X_2) + (X+Y),\adveq$$
so $(X_1, Y_2, X+Y)$ and $(X_1 - X_2, X+Y)$
each separately determine $(X_1 + Y_2, X+Y)$. Hence the submodularity inequality yields
$$ \Eta\{X_1, X_2, X+Y\} + \Eta\{X_1+Y_2, X+Y\} \le \Eta\{X_1, Y_2, X+Y\} + \Eta\{X_1- X_2, X+Y\}.\adveq$$
From $(X_1, X_2, X+Y)$
we can calculate $Y_1 = X+Y-X_1$ and $Y_2 = X+Y-X_2$, so this triple and the triple
$(X_1, X_2, Y_1, Y_2)$ determine each other. So the first term above is simply the additive energy
between $X$ and $Y$; that is
$$ \Eta\{X_1, X_2, X+Y\} = \Eta\{X_1, X_2, Y_1, Y_2\} = 2\Eta\{X,Y\} - \Eta\{X+Y\}.\adveq$$
Now since $X_1$ and $Y_2$ are conditionally independent relative to $X+Y$, we have
$$\eqalign{
\Eta\{X_1, Y_2, X+Y\} &= \Eta\{X_1, X+Y\} + \Eta\{Y_2, X+Y\} - \Eta\{X+Y\} \cr
&= \Eta\{X, X+Y\} + \Eta\{Y, X+Y\} - \Eta\{X+Y\} \cr
&= 2\Eta\{X,Y\} - \Eta\{X+Y\} \cr
}\adveq$$
For the last term above we split
$$\Eta\{X_1 - X_2, X+Y\} = \Eta\{X_1 - X_2 \given X+Y\} - \Eta\{X+Y\}.\adveq$$
Putting everything together, we obtain
$$\eqalign{
2\Eta\{X,Y\} - \Eta\{X+&Y\} + \Eta\{X_1 + Y_2, X+Y\} \cr
&\le 2\Eta\{X,Y\} + \Eta\{X_1-X_2\given X+Y\} - 2\Eta\{X+Y\},\cr
}\adveq$$
so that
$$\Eta\{ X_1 + Y_2\given X+Y\} \le \Eta\{X_1 - X_2\given X+Y\} - 2\Eta\{X+Y\}\le \Eta\{X_1 - X_2\}.\adveq$$
The previous lemma then gives
$$\Eta\{ X_1 + Y_2\given X+Y\} \le 2\Eta\{X\} + 2\Eta\{Y\} - \AA\{X,Y\},\adveq$$
and from our lower bound on $\AA\{X,Y\}$, we conclude that
$$\Eta\{ X_1 + Y_2\given X+Y\} \le {1\over 2} \Eta\{X\} + {1\over 2} \Eta\{Y\} + \log C,\adveq$$
which is what we wanted to show.\slug

Comparing constants with~\bref{18}, if we had started with an error term of $-3\log C$ we
would end up with error terms of $-6\log C$ in the bounds on $\Eta\{X_1\given X+Y\}$
and $\Eta\{Y_2\given X+Y\}$, whereas the 2010 paper only gives errors of $-\log C$. In return however,
our final upper bound only has an error of $3\log C$ where the other paper incurs a
corresponding error of $7\log C$.

\medskip\boldlabel Symmetry and asymmetry.
We note that the symmetry in
the hypothesis of the entropic Balog--Szeme\-r\'edi--Gowers theorem forces
the entropies of $X$ and $Y$ not to differ too wildly.

\proclaim Proposition {\advthm}. If $X$ and $Y$ are such that
$$\AA\{X,Y\} \ge {3\over 2} \Eta\{X\} + {3\over 2}\Eta\{Y\} - \log C,$$
then
$$\Eta\{X\} - 2\log C \le \Eta\{Y\} \le \Eta\{X\} + 2\log C.\adveq$$

\proof Using the bound~\refeq{\eqindepmax} and multiplying through by $2$,
we have
$$0\le 2\Eta\{X+Y\} \le \Eta\{X\} + \Eta\{Y\} + 2\log C.\adveq$$
So $\Eta\{X\} \le \Eta\{Y\} + 2\log C$ and $\Eta\{Y\} \le \Eta\{X\} + 2\log C$
both hold.\slug

In the most asymmetric case, where, say, $\Eta\{X\}$ is much greater than $\Eta\{Y\}$,
we can instead bound the distance between $\Eta\{X+Y\}$ and $\Eta\{Y\}$.

\proclaim Proposition {\advthm}. If $X$ and $Y$ are such that
$$\AA\{X,Y\} \ge 2 \Eta\{X\} + \Eta\{Y\} - \log C,$$
then
$$\Eta\{X+Y\} - \Eta\{Y\} \le \log C.\adveq$$

\proof We have
$$2\Eta\{X\} + \Eta\{Y\} - \log C \le \AA\{X,Y\} \le 2\Eta\{X\} + 2\Eta\{Y\} - \Eta\{X+Y\},\adveq$$
and cancelling terms produces the desired inequality.\slug

\advsect The small-energy regime
\hldest{xyz}{}{small}

We now perform a similar analysis to that of the previous section,
but in the case that the entropic additive energy is small. For the sumset setting we
have $E(A,B)\ge |A|\cdot|B|$, and in the entropic setting we have
$$\AA\{X,Y\} = 2\Eta\{X,Y\} - \Eta\{X+Y\} \ge \Eta\{X,Y\},\adveq$$

For sumsets, small additive energy implies that $|A+B|$ is large, as this counterpart
to Proposition~{\propinversebalog} shows. It also follows easily from the Cauchy--Schwarz inequality.

\edef\propsmallenergy{\the\thmcount}
\proclaim Proposition \advthm.
Let $A$ and $B$ be finite subsets of an abelian group.
If
$$E(A,B)\le C|A|\cdot |B|$$
for some constant $C$, then
$$|A+B| \ge {|A|\cdot|B|\over C}.\noskipslug\adveq$$

Once again, the implication in the entropic setting goes the wrong way around.

\edef\propnaivetwo{\the\thmcount}
\proclaim Proposition \advthm. Let $X$ and $Y$ be discrete random variables taking values in the
same abelian group. If
\global\edef\eqsumboundtwo{\the\eqcount}
$$\Eta\{X+Y\} \ge \Eta\{X\} + \Eta\{Y\} - \log C,\adveq$$
for some constant $C$, then
\global\edef\eqenergyboundtwo{\the\eqcount}
$$\AA\{X,Y\} \le \Eta\{X\} + \Eta\{Y\} + \log C.\adveq$$
If one adds the further assumption that
$\Eta\{X,Y\} \ge \Eta\{X\} + \Eta\{Y\} - C$,
then~{\rm \refeq{\eqenergyboundtwo}} implies~{\rm \refeq{\eqsumboundtwo}} with a worse constant, namely,
we may only conclude
\global\edef\eqworseboundtwo{\the\eqcount}
$$\Eta\{X+Y\} \ge \Eta\{X\} + \Eta\{Y\} - \log C-2C',\adveq$$
In particular, if $X$ and $Y$ are independent, then we can recover~{\rm \refeq{\eqsumboundtwo}}
from~{\rm \refeq{\eqenergyboundtwo}}.

\proof Assuming~\refeq{\eqsumboundtwo}, we have
$$\eqalign{
\AA\{X,Y\} &= 2\Eta\{X,Y\} - \Eta\{X+Y\} \cr
&\le 2\Eta\{X,Y\} - \Eta\{X\} - \Eta\{Y\} - \log C \cr
&\le \Eta\{X\} + \Eta\{Y\} - \log C.\cr
}\adveq$$
On the other hand, assuming~\refeq{\eqenergyboundtwo} as well as the lower bound on $\Eta\{X,Y\}$, we can bound
$$\Eta\{X+Y\}= 2\Eta\{X,Y\} -\AA\{X,Y\} \ge \Eta\{X\} + \Eta\{Y\} - \log C -2C'.\noskipslug\adveq$$

This time, the bound on $\Eta\{X+Y\}$ is more powerful, since it immediately implies
that $\Eta\{X,Y\} \ge \Eta\{X\} + \Eta\{Y\} - \log C$, which asserts a weak form of independence,
whereas for the other direction we must add this as a supplementary assumption.

\medskip\boldlabel Doubling constants and Sidon sets.
So far we have discussed the sumset $A+B$ of two sets. We now turn to the
special case where $A=B$. In this setting one speaks of the {\it doubling constant} $\sigma(A) = |A+A|/|A|$
of $A$. Given a random variable $X$ taking values in an abelian group, the {\it entropic doubling constant}
$\ss\{X\}$ is the quantity $\Eta\{X+X'\} - \Eta\{X\}$, where $X'$ is an independent copy of $X$.
This definition (or precisely, the exponential thereof) was established and studied in~\bref{18}.

From the definition~\refeq{\eqenergyonesetdef} of $\AA\{X\}$, we see that
\edef\eqdoublingformula{\the\eqcount}
$$\AA\{X\} = 3\Eta\{X\} - \ss\{X\}, \adveq$$
and hence quantitative statements about $\ss\{X\}$ are statements about $\AA\{X\}$ and vice versa.
For instance, it is easy to derive the special cases of Propositions~{\propnaive} and~{\propnaivetwo}
that result from setting $Y = X'$; namely,
$$ \AA\{X\} \ge 3\Eta\{X\} - \log C\quad\hbox{if and only if}\quad\ss\{X\}\le \log C\adveq$$
and
$$ \ss\{X\} \ge \Eta\{X\} - \log C\quad\hbox{if and only if}\quad\AA\{X\} \le 2\Eta\{X\} + \log C.\adveq$$
Furthermore, combining~\refeq{\eqdoublingformula} with the upper bound of Proposition~{\propbasicbounds}
yields a simpler proof of the inequality
$$\ss\{U_A\} \ge 3\log |A| - E(A)\adveq$$
when $U_A$ is the uniform distribution on $A$. This inequality appears as~\refeq{1{\rm .}1} in~\bref{9},
and two proofs of it are given (one using the weighted AM-GM inequality and the other using the monotonicity
of R\'enyi entropy).

Let us now discuss the extremal cases. In the combinatorial setting, the doubling constant is minimum
for cosets of subgroups, and Theorem~1.11 of~\bref{18} asserts that the same thing holds
in the entropic setting ($\ss\{X\} = 0$ if and only if $X$ is uniform on a coset of a finite subgroup).

At the other end of the spectrum, the doubling constant is maximum for {\it Sidon sets}. These are
sets $A$ whose pairwise sums are all distinct; that is, $a+b = a'+b'$ implies that $\{a,b\} = \{a',b'\}$.
Returning to entropies, we have the equivalent bounds $\AA\{X\} \ge 2\Eta\{X\}$
and $\ss\{X\} \le \Eta\{X\}$.
These bounds are satisfied with equality when $X$ takes only one value. Furthermore, Sidon sets
give rise to random variables $X$ that come close to attaining the bound.

\edef\propsidon{\the\thmcount}
\proclaim Proposition {\advthm}. Let $X$ be any (discrete) distribution on a Sidon set $A$ in an abelian
group. Then $\ss\{X\} \ge \Eta\{X\} - 1$ (hence $\AA\{X\} \le 2\Eta\{X\} - 1$).

\proof Let $X'$ be an independent copy of $X$. Since $A$ is Sidon, if we know that $X+X' = a$,
then we know the values of $X$ and $X'$ up to ordering. Impose a total order $\le$ on $A$ and let
$$ Y = \cases{ 0, & if $X\le X'$;\cr 1, & otherwise.}$$
Now the pairs $(X,X')$ and $(Y, X+X')$ determine each other, so
$$ 2\Eta\{X\} = \Eta\{X,X'\} = \Eta\{Y, X+X'\} \le \Eta\{Y\} + \Eta\{X+X'\} \le 1 + \Eta\{X+X'\},\adveq$$
where in the rightmost inequality we have used the fact that $Y$ takes values in a $2$-element set.
The result then follows upon rearranging.\slug

In light of this result, let us say that a random variable $X$ taking values in an abelian group
is {\it Sidon}
if $\ss\{X\} \ge \Eta\{X\} - 1$. This is strictly more general than requiring that the range of $X$
be a Sidon set, as shown by the following small example. Let $X$ and $X'$ be independently
uniform on the set $\{0,1,2\}\subseteq \ZZ$. We have $\Eta\{X\} = \log 3$, and using the formula
$$\Eta\{Z\} = \sum_z \pr\{Z=z\} \log \biggl({1\over\pr\{Z=z\}}\biggr),\adveq$$
we compute
$$\Eta\{X+X'\} = {2\over 9} \log 9 + {4\over 9}\log(9/2) + {1\over 3}\log 3 - \log 3 > 0.4244.\adveq$$
Since $\Eta\{X\} - 1 = \log 3 - 1 < 0.098$, $X$ is a Sidon random variable, but $\{0,1,2\}$ is not a Sidon set
because $0+2 = 1+1$.

It is easy to verify that if $A$ and $B$ are sets with $B\subseteq A$ and $B$ is Sidon, then
$|B| \le \sqrt{2|A+A|}$. We prove the entropic version of this statement.

\edef\propsidonconditioning{\the\thmcount}
\proclaim Proposition {\advthm}. Let $X$ be a random variable supported on an abelian group, and let $Z$
be any random variable such that $\Eta\{X\given Z\}$ is Sidon. Then
$$\Eta\{X\given Z\} \le {\Eta\{X+X'\} + 1\over 2}.\adveq$$

\proof By definition of a Sidon random variable,
$$\Eta\{X+X'\given Z\} - \Eta\{X\given Z\} \ge \Eta\{X\given Z\} - 1,\adveq$$
so
$$2\Eta\{X\given Z\} \le \Eta\{X+X'\given Z\} + 1 \le \Eta\{X+X'\} + 1,\adveq$$
and dividing through by $2$ finishes the proof.\slug

It should be mentioned that our investigation of random variables with large entropic doubling
was greatly expanded upon in Section~5 of~\bref{13}.

\advsect Multiplicative energy and the sum-product problem
\hldest{xyz}{}{mult}

The setup that we used to construct the entropic additive energy can be used to define a multiplicative
energy as well, in the case that the random variables under scrutiny take values in a ring.
Concretely, given $X$ and $Y$ taking values in a ring,
we let $(X_1,Y_1)$ and $(X_2, Y_2)$ be conditionally independent trials
of $(X,Y)$ relative to $XY$, and define the {\it entropic multiplicative energy} to be
$$\MM\{X,Y\} = \Eta\{X_1, Y_1, X_2, Y_2\}.$$
As before, this reduces to the simpler formula
$$\MM\{X,Y\} = 2\Eta\{X,Y\} - \Eta\{X Y\}.$$
This definition does not require the existence of inverses nor commutativity.
but in the case that the ring is
a (commutative) field $K$, then by letting the group be $K^*$ under multiplication, all the results we proved
for $\AA\{X,Y\}$ in general groups carry through for $\MM\{X,Y\}$, so long as $X$ and $Y$ take values in $K^*$.

One might wonder if there is a sum-product phenomenon for entropies of random variables (see, e.g, Chapter~2.8
of~\bref{17} for the corresponding additive-combinatorial theory). That is, if $X$
has large multiplicative energy $\MM\{X\}$, must its additive energy $\AA\{X\}$ be small?
In general, this is not true.
Let $K$ be a finite field, let $K'$ be a subfield of $K$ of cardinality $q$,
and let $X$ and $X'$ be independent copies of $U_{{K'}^*}$.
Then $\Eta\{X X'\} = \Eta\{X\} = \log (q-1)$, so $\MM\{X\} = 3\Eta\{X\}$ is large. On the other hand,
letting $A$ be the support of $X+X'$, we have
$$\eqalign{
\Eta\{X+ X'\} &=\sum_{x\in A} \pr\{X+  X' =x\} \log \biggl({1\over \pr\{X+ X' = x\}}\biggr) \cr
&= {1\over q-1} \log (q-1) + (q-1) {q-2\over (q-1)^2} \log\biggl({(q-1)^2\over q-2}\biggr) \cr
&= 2 \log (q-1) - \log (q-2) + O\biggl({\log q\over q}\biggr),\cr
}\adveq$$
so $\AA\{X\} = 2\log (q-1) + \log (q-2) - O\bigl(( \log q)/q\bigr) \approx 3\Eta\{X\}$ is also large.

In the remainder of this paper, we formulate some conjectures.
The field $\FF_p$ has no nontrivial subfields, so it is easy to forbid the scenario above using a quantitative
stipulation on $\Eta\{X\}$.
Once this is done, we conjecture that either $\AA\{X\}$ or $\MM\{X\}$ must be small. This is by analogy
with the following theorem
of J.~Bourgain, N.~Katz, and T.~Tao~\bref{3}.

\edef\thmbourgainkatztao{\the\thmcount}
\parenproclaim Theorem {\advthm} (Bourgain--Katz--Tao, {\rm 2004}).
Let $A\subseteq \FF_p^*$ for some prime $p$, and suppose that
$$p^\delta < |A| < p^{1-\delta}$$
for some $\delta>0$. Then there exists $\eps$ depending only on $\delta$ such that
$$\max\bigl(|A+A|, |A\cdot A|\bigr) > |A|^{1+\eps}.\adveq$$

Our conjectured entropic analogue has the following form.

\edef\conjbourgainkatztao{\the\thmcount}
\parenproclaim Conjecture {\advthm} (Entropic Bourgain--Katz--Tao theorem in $\FF_p$).
Let $0< \delta<1$.
There exist constants $\eps$ and $p_0$ depending only on $\delta$ such that
if $p$ is a prime greater than $p_0$,
then any random variable $X$ taking values in $\FF_p^*$ with
$$\delta \log p \le \Eta\{X\} \le (1-\delta)\log p$$
and
$$\MM\{X\}\ge (3-\eps)\Eta\{X\}$$
must satisfy
$$\AA\{X\} < (3-\eps)\Eta\{X\}.$$

First we show that our conjecture implies the original Bourgain--Katz--Tao theorem for large enough $p$
(depending on $\delta$).

\medskip\noindent{\it Proof that Conjecture~{\conjbourgainkatztao} implies Theorem~{\thmbourgainkatztao} for large
$p$}.\enspace
Let $0<\delta<1$, let $p\ge p_0$, and let
$A$ be a subset of $\FF_p^*$ such that $p^\delta < |A| < p^{1-\delta}$.
Then letting $U_A$ be the indicator random variable on $A$, we have $\Eta\{U_A\} = \log |A|$, so
$$\delta \log p \le \Eta\{U_A\} \le (1-\delta)\log p\adveq$$
By Conjecture~{\conjbourgainkatztao}, either $\MM\{U_A\}$ or $\AA\{U_A\}$ must be less than $(3-\eps)\Eta\{U_A\}$,
so either $\Eta\{U_A + U_A'\}$ or $\Eta\{U_A U_A'\}$ is greater than $(1+\eps)\Eta\{U_A\}$,
where $U_A'$ is another indicator random variable on $A$, independent of $U_A$.
Suppose that the sum $U_A+U_A'$ has high entropy. Then
$$\log |A+A| \ge \Eta\{U_A+U_A'\} > (1+\eps)\Eta\{U_A\} = (1+\eps)\log |A|,\adveq$$
and consequently
$$|A+A|> |A|^{1+\eps}.\adveq$$
If instead the product $U_A U_A'$
has high entropy then {\it mutatis mutandis}
we conclude that $|A\cdot A|> |A|^{1+\eps}$.\slug

We were unable to prove the entropic version of the Bourgain--Katz--Tao theorem,
but can show that it follows from a related conjecture,
namely, an entropic version of a lemma attributed to A.~Glibichuk and S.~V.~Konyagin
(see, e.g.,~\bref{8}).

\edef\conjkonyagin{\the\thmcount}
\parenproclaim Conjecture {\advthm} (Entropic Glibichuk--Konyagin lemma).
There is an integer $k\ge 3$ such that the following holds.
Let $X$ be a random variable taking values in $\FF_p$ for a prime $p$. Then for independent copies
$X_1, \ldots, X_k, X_1', \ldots, X_k'$ of $X$, we have
$$\Eta\{X_1 X_1' + \cdots + X_k X_k'\} \ge \min\bigl(2\Eta\{X\}, \log p\bigr) - 1.\adveq$$

In the original additive-combinatorial version of this statement, $k$ can be taken to be $3$, but any
fixed integer is enough for our purposes.

We begin by observing that the following
entropic analogue of the Katz--Tao lemma~\bref{12} follows
from a recent theorem of A.~Math\'e and W.~L.~O'Regan~\bref{14}.
The first claim in the following lemma
was stated as Theorem~4.3 in a preliminary version of~\bref{14} before being removed in
the journal version of that paper. For completeness, we supply a full proof,
while also isolating the particular case that is of interest to us.

\edef\lemkatztao{\the\thmcount}
\parenproclaim Lemma {\advthm} (Entropic Katz--Tao lemma).
Let $X$, $Y$, $Z$, and $W$ be independent random variables taking values in $K^*$ for some
commutative ring $K$. Then
\global\edef\eqmoclaim{\the\eqcount}
$$\eqalign{
\Eta\{XY + Z&W\} + \Eta\{X\} + \Eta\{Y\} + 2\Eta\{Z\} + 2\Eta\{W\} \cr
&\le \Eta\{X+Y\} + \Eta\{Z-W\} + \Eta\{XZ\} + 2\Eta\{YZ\} + 2\Eta\{ZW\}.\cr
}\adveq$$
In particular, if $X$, $Y$, $Z$, and $W$ are identically distributed, with
$$\min\bigl(\AA\{X\},\MM\{X\}\bigr) \ge 3\Eta\{X\} - \log C$$
for some parameter $C>0$, then
\global\edef\eqmoclaimtwo{\the\eqcount}
$$\Eta\{XY + ZW\} \le \Eta\{X\} + 8\log C.\adveq$$

\proof Let $S = \bigl( XY, ZW, (X+Z)Y, Z(Y-W)\bigr)$ for short. The tuples $(S,Y)$ and
$(X,Y,Z,W)$ determine each other, so $\Eta\{S,Y\} = \Eta\{X,Y,Z,W\}$. By submodularity,
$$\Eta\{S,Y\} + \Eta\{XY\} \le \Eta\{S\} + \Eta\{X,Y\},$$
so
\edef\eqfirstsub{\the\eqcount}
$$\eqalign{
\Eta\{S\} &\ge \Eta\{X,Y,Z,W\} + \Eta\{XY\} - \Eta\{X,Y\} \cr
&= \Eta\{Z,W\} + \Eta\{XY\}, \cr
}\adveq$$
by independence.
Another application of submodularity yields
$$\Eta\{XY+ZW\} + \Eta\{S\} \le \Eta\{XY, ZW\} + \Eta\bigl\{ (X+Z)Y, Z(Y-W)\bigr\},$$
since $(X+Z)Y + Z(Y-W) = XY+ZW$. Substituting in~\refeq{\eqfirstsub}, we get
\edef\eqsecondsub{\the\eqcount}
$$\Eta\{XY+ZW\} + \Eta\{Z,W\} \le \Eta\{ZW\} + \Eta\bigl\{ (X+Z)Y\bigr\} + \Eta\bigl\{Z(Y-W)\bigr\},\adveq$$
where we have cancelled a $\Eta\{XY\}$ term from both sides.
Lastly, note that
$$\eqalign{
\Eta\bigl\{ (X+Z)Y\bigr\} &\le \Eta\{XZ, YZ\} + \Eta\{X+Z,Y\} - \Eta\{X,Y,Z\} \cr
&\le \Eta\{XZ\} + \Eta\{YZ\} + \Eta\{X+Z\} - \Eta\{X\} - \Eta\{Z\} \cr
}$$
by submodularity, and likewise
$$\eqalign{
\Eta\bigl\{ Z(Y-W)\bigr\} &\le \Eta\{YZ, ZW\} + \Eta\{Z,Y-W\} - \Eta\{Y,Z,W\} \cr
&\le \Eta\{YZ\} + \Eta\{ZW\} + \Eta\{Y-W\} - \Eta\{Y\} - \Eta\{W\}.\cr
}$$
Substituting both of these inequalities into~\refeq{\eqsecondsub} gives exactly~\refeq{\eqmoclaim}.

The hypothesis in the second claim implies that
$$\max\bigl(\Eta\{X+Y\}, \Eta\{XY\}\bigr) \le \Eta\{X\} + \log C.\adveq$$
But by~\refeq{\eqmoclaim} one has
$$\Eta\{XY+ZW\} \le 5\Eta\{XY\} + 3\Eta\{X+Y\} - 7\Eta\{X\},\adveq$$
and~\refeq{\eqmoclaimtwo} follows easily.\slug

Assuming the conjectured entropy version of the Glibichuk--Konyagin result,
the rest of the proof outlined in~\bref{8}
carries over from the classical setting to the entropic one without much ado.
In the proof, we will invoke the entropy version
of the Pl\"unnecke--Ruzsa inequality, which is due to V.~A.~Kaimanovich and
A.~M.~Vershik~\bref{11}.

\edef\thmplunnecke{\the\thmcount}
\parenproclaim Theorem {\advthm} (Entropic Pl\"unnecke--Ruzsa inequality).
Let $W, W_1, \ldots, W_m$ be independent random variables with finite entropy taking values in the same
abelian group. If
$$\Eta\{W+W_i\} \le \Eta\{W\} + \log C_i$$
for $1\le i\le m$ and constants $C_1,\ldots, C_m$, then
$$\Eta\{W+W_1+\cdots+W_m\} \le \Eta\{W\} + \log(C_1\cdots C_m).\noskipslug\adveq$$

\medskip\noindent{\it Proof that Conjecture~{\conjkonyagin}
implies Conjecture~{\conjbourgainkatztao}}.\enspace
Let $0<\delta<1$ and let $\eps = \eps(\delta)$ be chosen later.
Suppose, for a contradiction, that $\AA\{X\}$ and $\MM\{X\}$ are both at least $(3-\eps)\Eta\{X\}$.
By Lemma~{\lemkatztao} with $\log C = \eps\Eta\{X\}$, letting
$X_1$, $X_2$, $X_3$, and $X_4$ are independent copies of $X$ and setting $Y = X_1 X_2 + X_3 X_4$
for brevity, we have
\edef\eqyupper{\the\eqcount}
$$\Eta\{Y\} \le (1+8\eps)\Eta\{X\}.\adveq$$
Next we work towards a lower bound on $\Eta\{Y\}$.
Now we apply Conjecture~{\conjkonyagin}. Let $k$ be the integer from that statement and let
$$Y_k = \bar X_1 \widetilde X_1 + \cdots + \bar X_k  \widetilde X_k\adveq$$
for independent copies $\bar X_1,\ldots,\bar X_k,\widetilde X_1,\ldots,\widetilde X_k$ of $X$.
If
$${1\over 2} \log p \le \Eta\{X\} \le (1-\delta)\log p,\adveq$$
then setting $p_0 = 2^{4/\delta}$, we have
$${\delta\over 2} \Eta\{X\} \ge {\delta\over 4} \log p \ge 1\adveq$$
for all $p\ge p_0$. Assuming for the rest of the proof that $p$ satisfies this bound, we have
$$\eqalign{
\Eta\{Y_k\} &\ge \log p - 1 \cr
&\ge {1\over 1-\delta} \Eta\{X\} - 1 \cr
&\ge (1+\delta)\Eta\{X\} -1 \cr
&\ge \biggl(1+{\delta\over 2}\biggr) \Eta\{X\}.\cr
}\adveq$$
On the other hand, if $\Eta\{X\} \le (1/2)\log p$, then we simply use the bound
$$\Eta\{X\} \ge \delta \log p \ge 4$$
to get
$$\Eta\{Y_k\} \ge 2\Eta\{X\} - 1\ge {7\over 4} \Eta\{X\}
\ge \biggl(1+{\delta\over 2}\biggr) \Eta\{X\},\adveq$$
since $\delta < 1$. Thus in either case we have
$\Eta\{Y_k\} \ge (1+\delta/2)\Eta\{X\}$.
Now $\MM\{X\} \ge (3-\eps)\Eta\{X\}$ gives $\Eta\{X_1 X_2\} \le (1+\eps)\Eta\{X\}$, so
\edef\eqdeltathree{\the\eqcount}
$$\Eta\{Y_k\} \ge {1+\delta/2\over 1+\eps}\Eta\{X_1 X_2\}.\adveq$$
Let $\eps < \delta/\bigl(18(k-1)\bigr)$, so that
$$\eqalign{
\Eta\{Y_k\} &\ge {1+\delta/2\over 1+\delta/\bigl(18(k-1)\bigr)}\Eta\{X_1 X_2\} \cr
&\ge {1+\delta/2\over 1+\delta/36}\Eta\{X_1 X_2\} \cr
&\ge {1+4\delta/9 + \delta/36 + \delta^2/81 \over 1+\delta/36}\Eta\{X_1 X_2\} \cr
&= \biggl(1+{4\delta\over 9}\biggr)\Eta\{X_1 X_2\}.\cr
}\adveq$$
But letting $L$ satisfy
$$\Eta\{Y\} = \Eta\{X_1 X_2\} + \log L,$$
by the entropic Pl\"unnecke--Ruzsa inequality with $m=k-1$, the variables $W, W_1,\ldots, W_{k-1}$
set to independent copies of $X_1 X_2$, and $C_i = L$ for all $1\le i <k$,
we obtain
$$\Eta\{Y_k\} \le \Eta\{X_1 X_2\} + (k-1)\log L.\adveq$$
Combining this bound with~\refeq{\eqdeltathree} yields
$$\log L \ge {4\delta\over 9(k-1)} \Eta\{X X'\}.\adveq$$
This means that
$$\Eta\{Y\} \ge \biggl(1+{4\delta\over 9(k-1)}\biggr)\Eta\{X_1 X_2\},\adveq$$
and chaining this bound with~\refeq{\eqyupper} gives us
$$(1+8\eps)\Eta\{X\} \ge \biggl(1+{4\delta\over 9(k-1)}\biggr) \Eta\{X\},\adveq$$
contradicting our earlier assumption that $\eps < \delta/\bigl(18(k-1)\bigr)$.\slug

On the real line, it was conjectured by P.~Erd\H{o}s and E.~Szemer\'edi~\bref{4}
that any finite set $A\subseteq\RR$ has
$\max\bigl( |A+A|,|A\cdot A|\bigr) \ge |A|^{2-o(1)}$.
The most recent lower bound, $|A|^{4/3 + 2/951 - o(1)}$, was established in 2025 by
T.~Bloom~\bref{2}, improving slightly upon the exponent of
$4/3 + 2/1167$ obtained by M.~Rudnev and S.~Stevens~\bref{15} in 2020.
We make the following tentative conjecture in the entropy setting.

\edef\conjreal{\the\thmcount}
\proclaim Conjecture~{\advthm}. There exists $\eps>0$ such that for any $X$ taking values in $\RR$ with
$\Eta\{X\} < \infty$,
$$\max\bigl( \AA\{X\}, \MM\{X\}\bigr) \le \bigl(3-\eps+o(1)\bigr)\Eta\{X\}.\adveq$$

As we saw earlier over $\FF_p$, the truth of this conjecture for some $\eps>0$ immediately implies the
analogous bound in the classical setting.
Following the initial preprint release of the present paper, however,
it was shown in~\bref{13}
that if Conjecture~{\conjreal} is true for some $\eps>0$, one must have $\eps\le 1/3$. In other words,
it is not possible to improve on Bloom's exponent in the classical Erd\H{o}s--Szemer\'edi problem by
means of this entropy analogue.
Nonetheless we leave our conjecture up in the hope that it, along
with other entropic sum-product problems, should garner independent interest.

\section Acknowledgements
\hldest{xyz}{}{acks}

The author wishes to thank Ben Green for alerting his attention to the cited paper of M\'ath\'e and O'Regan,
and Hamed Hatami for numerous helpful comments on a preliminary version of this manuscript.
Thanks are also due to the two anonymous referees, whose comments and suggestions have been used to
greatly improve the paper.
The author is funded by the Natural Sciences and Engineering Research Council of Canada.

\section References
\hldest{xyz}{}{refs}

\parskip=0pt
\hyphenpenalty=-1000 \pretolerance=-1 \tolerance=1000
\doublehyphendemerits=-100000 \finalhyphendemerits=-100000
\frenchspacing
\def\bref#1{[#1]}
\def\beginref{\noindent}
\def\endref{\medskip}
\vskip\parskip

\beginref
\parindent=20pt\item{\bref{1}}
\hldest{xyz}{}{bib1}%
Antal Balog
and Endre Szemer\'edi,
``A statistical theorem of set addition,''
{\sl Combinatorica}\/
{\bf 14}
(1994),
263--268.
\endref
\beginref
\parindent=20pt\item{\bref{2}}
\hldest{xyz}{}{bib2}%
Thomas Bloom,
``Control and its applications in additive combinatorics,''
{\sl arXiv:2501.09470}\/
(2025),
28~pp.
\endref
\beginref
\parindent=20pt\item{\bref{3}}
\hldest{xyz}{}{bib3}%
Jean Bourgain,
Nets Katz,
and Terence Tao,
``A sum-product estimate in finite fields, and applications,''
{\sl Geometric and Functional Analysis}\/
{\bf 14}
(2004),
27--57.
\endref
\beginref
\parindent=20pt\item{\bref{4}}
\hldest{xyz}{}{bib4}%
Paul Erd\H{o}s
and Endre Szemer\'edi,
``On sums and products of integers,''
{\sl Studies in Pure Mathematics}\/
(1983),
213--218.
\endref
\beginref
\parindent=20pt\item{\bref{5}}
\hldest{xyz}{}{bib5}%
William Timothy Gowers,
``A new proof of Szemer\'edi's theorem for arithmetic progressions of length four,''
{\sl Geometric and Functional Analysis}\/
{\bf 8}
(1998),
529--551.
\endref
\beginref
\parindent=20pt\item{\bref{6}}
\hldest{xyz}{}{bib6}%
William Timothy Gowers,
Ben Green,
Freddie Manners,
and Terence Tao,
``On a conjecture of Marton,''
{\sl Annals of Mathematics}\/
{\bf 201}
(2025),
515--549.
\endref
\beginref
\parindent=20pt\item{\bref{7}}
\hldest{xyz}{}{bib7}%
William Timothy Gowers,
Ben Green,
Freddie Manners,
and Terence Tao,
``Marton's conjecture in abelian groups with bounded torsion,''
{\sl Annales de la Facult\'e des Sciences de Toulouse}\/
{\bf 2026}
(2026),
1--33.
\endref
\beginref
\parindent=20pt\item{\bref{8}}
\hldest{xyz}{}{bib8}%
Ben Green,
``Sum-product phenomena in ${\bf F}_p$: a brief introduction,''
{\sl arXiv: 0904.2075}\/
(2009),
10~pp.
\endref
\beginref
\parindent=20pt\item{\bref{9}}
\hldest{xyz}{}{bib9}%
Ben Green,
Freddie Manners,
and Terence Tao,
``Sumsets and entropy revisited,''
{\sl Random Structures and Algorithms}\/
{\bf 66}
(2025),
e21252.
\endref
\beginref
\parindent=20pt\item{\bref{10}}
\hldest{xyz}{}{bib10}%
Peter Hegarty,
``Some explicit constructions of sets with more sums than differences,''
{\sl Acta Arithmetica}\/
{\bf 130}
(2007),
61--77.
\endref
\beginref
\parindent=20pt\item{\bref{11}}
\hldest{xyz}{}{bib11}%
Vadim Adolfovich Kaimanovich
and Anatoly Moiseevich Vershik,
``Random walks on discrete groups: boundary and entropy,''
{\sl Annals of Probability}\/
{\bf 11}
(1983),
457--490.
\endref
\beginref
\parindent=20pt\item{\bref{12}}
\hldest{xyz}{}{bib12}%
Nets Katz
and Terence Tao,
``Some connections between the Falconer and Furstenburg conjectures,''
{\sl New York Journal of Mathematics}\/
{\bf 7}
(2001),
148--187.
\endref
\beginref
\parindent=20pt\item{\bref{13}}
\hldest{xyz}{}{bib13}%
Rupert Li,
Lampros Gavalakis,
and Ioannis Kontoyiannis,
``Entropic additive energy and entropy inequalities for sums and products,''
{\sl {\rm To appear in} IEEE Transactions on Information Theory}\/
(2025),
26~pp.
\endref
\beginref
\parindent=20pt\item{\bref{14}}
\hldest{xyz}{}{bib14}%
Andr\'as M\'ath\'e
and William Lewis O'Regan,
``Discretised sum-product theorems by Shannon-type inequalities,''
{\sl Journal of the London Mathematical Society}\/
{\bf 112}
(2025),
e70389.
\endref
\beginref
\parindent=20pt\item{\bref{15}}
\hldest{xyz}{}{bib15}%
Misha Rudnev
and Sophie Stevens,
``An update on the sum-product problem,''
{\sl Mathematical Proceedings of the Cambridge Philosophical Society}\/
{\bf 173}
(2020),
411--430.
\endref
\beginref
\parindent=20pt\item{\bref{16}}
\hldest{xyz}{}{bib16}%
Imre Zolt\'an Ruzsa,
``Sumsets and entropy,''
{\sl Random Structures and Algorithms}\/
{\bf 34}
(2009),
1--10.
\endref
\beginref
\parindent=20pt\item{\bref{17}}
\hldest{xyz}{}{bib17}%
Terence Tao
and Van Ha Vu,
{\sl Additive Combinatorics}
(Cambridge:
Cambridge University Press,
2006).
\endref
\beginref
\parindent=20pt\item{\bref{18}}
\hldest{xyz}{}{bib18}%
Terence Tao,
``Sumset and inverse sumset theory for Shannon entropy,''
{\sl Combinatorics, Probability, and Computing}\/
{\bf 19}
(2010),
603--639.
\endref
\beginref
\parindent=20pt\item{\bref{19}}
\hldest{xyz}{}{bib19}%
Andr\'e Weil,
``De la m\'etaphysique aux math\'ematiques,''
{\sl Sciences}\/
{\bf 60}
(1960),
52--56.
\endref
\beginref
\goodbreak\bye